\documentclass[12pt,a4paper,draft]{article}
\usepackage[english]{babel}
\usepackage{srcltx}
\usepackage{cite}
\usepackage{latexsym,amsfonts,amssymb,amsmath,longtable,stmaryrd,mathrsfs,amsthm}
\usepackage[mathlines,pagewise]{lineno}
\renewcommand{\leq}{\leqslant}
\renewcommand{\geq}{\geqslant}
\renewcommand{\unlhd}{\trianglelefteqslant}

\renewcommand{\mathcal}{\mathscr}
\theoremstyle{plain}
\newtheorem{Lemma}{{\bfseries Lemma}}%}
\newtheorem{Cor}{{\bfseries Corollary}}
\newtheorem{Theo}{{\bfseries Theorem}}
\newtheorem{Prop}{{\bfseries Proposition}}
\theoremstyle{definition}
\newtheorem{Remark}{{\bfseries Remark}}
\newtheorem{Prob}{{\bfseries Problem}}
\newtheorem{Conj}{{\bfseries Conjecture}}

 \DeclareMathOperator{\Aut}{Aut} 
\DeclareMathOperator{\SL}{SL}

\DeclareMathOperator{\PSL}{PSL}
\DeclareMathOperator{\Sym}{Sym}

\modulolinenumbers[5]
\title{\vspace{-1cm} %\hfill{\normalsize MSC2010 20D20}{
Number of Sylow subgroups in finite groups\thanks{The first  author is supported by
a NNSF grant of China (Grant \# 11371335)
and Wu Wen-Tsun Key Laboratory of Mathematics of Chinese Academy of Sciences. The second
 author is supported by Chinese Academy of Sciences President's International Fellowship Initiative (Grant \#2017VMA0049)}}
\begin{document}
%\linenumbers

\author{Wenbin Guo\\
{\small Department of Mathematics, University of Science and
Technology of China,}\\ {\small Hefei 230026, P. R. China}\\
{\small E-mail:
wbguo@ustc.edu.cn}\\ \\
Evgeny P. Vdovin\\
{\small Sobolev Institute of Mathematics and Novosibirsk State University,}\\
{\small Novosibirsk 630090, Russia}\\
{\small E-mail: vdovin@math.nsc.ru}}

 \date{}
\maketitle

\pagenumbering{arabic}
\begin{abstract}

Denote by $\nu_p(G)$ the number of Sylow $p$-subgroups of $G$. It is not difficult to see that $\nu_p(H)\leq\nu_p(G)$
for $H\leq G$,  however $\nu_p(H)$ does not divide $\nu_p(G)$ in general. In this paper we reduce the question
whether $\nu_p(H)$ divides $\nu_p(G)$ for every $H\leq G$ to almost simple groups. This result substantially
generalizes the previous result by G. Navarro and also provides an alternative proof for the  Navarro theorem.


\end{abstract}

\section*{Introduction}

Throughout this paper, all groups are finite. Let $G$ be  a group and $p$ be a prime. Denote by $\nu_p(G)$ the number
of Sylow $p$-subgroups of $G$. It is a trivial exercise to check that $\nu_p(H)\leq \nu_p(G)$ for every subgroup $H$ of
$G$. However $\nu_p(H)$ does  not necessarily divide $\nu_p(G)$ in general. For example,  let $G=A_5$ and $H=A_4$, then
$\nu_3(G)=10$ and $\nu_3(H)=4$. In 2003, G. Navarro in \cite{Navarro} proved that if $G$ is $p$-solvable then
$\nu_p(H)$ divides $\nu_p(G)$ for every $H\leq G$.

We say that a group $G$ satisfies $\mathbf{DivSyl}(p)$ if $\nu_p(H)$ divides $\nu_p(G)$ for every $H\leq G$. Our goal
is to prove that $G$ satisfies $\mathbf{DivSyl}(p)$ provided every nonabelian composition factor of $G$ satisfies
$\mathbf{DivSyl}(p)$ (see Theorem~\ref{main} below). This result substantially generalizes the result by Navarro since
in $p$-solvable groups all nonabelian composition factors evidently satisfy $\mathbf{DivSyl}(p)$. Our technique can
also be applied to derive the Navarro theorem, so we provide an alternative proof for the Navarro theorem.

In order to formulate the main theorem, we need to recall some definitions.

Let $A,B,H$ be subgroups of  $G$ such that $B\unlhd A$. Define the {\em
normalizer} $N_H(A/B)$  of $A/B$ in $H$  by $N_H(A)\cap N_H(B)$. If $x\in
N_H(A/B)$, then $x$ induces an automorphism of $A/B$ acting by $Ba\mapsto B
x^{-1}ax$. Thus there exists a homomorphism  $N_H(A/B)\rightarrow \Aut(A/B)$.
The image of the homomorphism is denoted by $\Aut_H(A/B)$ and is called the
{\em group of $H$-induced automorphisms} on $A/B$, while the kernel of the
homomorphism is denoted by $C_H(A/B)$. If $B=1$, then we simply write
$\Aut_H(A)$. The groups of induced automorphisms were introduced by F. Gross
in \cite{GrossExistence}, where the author says that this notion is taken
from unpublished lectures by H. Wielandt. Clearly $C_H(A/B)=C_G(A/B)\cap H$,
so
$$\Aut_H(A/B)=N_H(A/B)/C_H(A/B)\simeq N_H(A/B)C_G(A/B)/C_G(A/B)\leq \Aut_G(A/B),$$  that is,
$\Aut_H(A/B)$ can be embedded into  $\Aut_G(A/B)$ in a natural way. Hence, without lose of generality, we may  assume that $\Aut_H(A/B)$ is a subgroup of
$\Aut_G(A/B)$.

A composition series is called an {\em $(rc)$-series}\footnote{This notion was introduced by V. A. Vedernikov in
\cite{Ved}} if it is a refinement of a chief series.

\begin{Theo}\label{main}
{\em (Main Theorem)} Let $$1=G_0<G_1<\ldots<G_n=G$$ be an $(rc)$-series of $G$. Assume that, for each nonabelian
$G_i/G_{i-1}$ and for every $p$-subgroup $P$ of $\Aut_G(G_i/G_{i-1})$, the group $P(G_i/G_{i-1})$ satisfies {\em
$\mathbf{DivSyl}(p)$}. Then $G$ satisfies  $\mathbf{DivSyl}(p)$.
\end{Theo}

At the end of the paper, we discuss possible way of improving this statement and also explain how Navarro's result can
be derived from this theorem.

\section{Preliminaries}

We write $H\leq G$ and $H\unlhd G$ if $H$ is a subgroup of $G$ and $H$ is a
normal subgroup of $G$ respectively, $p$ denotes a prime number. A natural
number $n$ is called a {\em $p$-number} (respectively a {\em $p'$-number}) if
it is a power of $p$ (respectively if it is coprime to $p$). We denote by
$n_p$ and $n_{p'}$ the $p$- and $p'$- part of the natural number $n$. A group
is said to be a {\em $p$-group} (respectively a {\em $p'$-group}) if its
order is a $p$-number (respectively a $p'$-number). A group $G$ is called
{\em $p$-solvable} if it possesses a subnormal series with all sections
either $p$- or $p'$-groups. If $\varphi:X\rightarrow Y$ is a map, then
$\varphi(x)$ denotes the image of $x\in X$, while for $\sigma\in\Sym_k$ we
write $i\sigma$ for the image of $i\in\{1,\ldots,k\}$ under $\sigma$. By
$Syl_p(G)$ we denote the set of all Sylow $p$-subgroups of~$G$. All
unexplained notations and definitions can be found in~\cite{Isaacs}.

\begin{Lemma}\label{nupfactor}
Assume that $A\unlhd G$. Then $\nu_p(G)=\nu_p(G/A)\cdot \nu_p(PA)$  for some (hence for every) $P\in Syl_p(G)$, and
$\nu_p(PA)=\vert A:N_A(P)\vert$.
\end{Lemma}

\begin{proof}
Consider a natural homomorphism
$$\overline{\phantom{G}}:G\rightarrow \overline{G}=G/A.$$
Clearly $\nu_p(PA)$ is the number of all  Sylow $p$-subgroups $Q$ of $G$ such
that image of $Q$ equals $\overline{P}$ since $\overline{P}=\overline{Q}$ if
and only if $Q\leq PA$. By Sylow theorem, for every $Q\in Syl_p(G)$, there
exists $x\in G$ such that $P^x=Q$, so
\begin{equation*}
(PA)^x=P^xA^x=P^xA=QA.
\end{equation*}
This shows that $\nu_p(PA)$ does not depend on the choice of $P$, and so it
is the same for any Sylow $p$-subgroup of $G$. Thus we obtain
$\nu_p(G)=\nu_p(G/A)\cdot \nu_p(PA)$. By Sylow theorem,
\begin{equation*}
\nu_p(PA)=\vert PA:N_{PA}(P)\vert =\vert A:N_A(P)\vert,
\end{equation*}
where the second identity follows from the Dedekind theorem.
\end{proof}

\begin{Lemma}\label{extension}
Assume that $A\unlhd G$ and  both $AP$ and $G/A$ satisfy $\mathbf{DivSyl}(p)$ for some (hence for every) $P\in
Syl_p(G)$. Then $G$ satisfies $\mathbf{DivSyl}(p)$.
\end{Lemma}

\begin{proof}
Let $H\leq G$. Choose $Q\in Syl_p(G)$ so that $Q\cap H\in Syl_p(H)$. Since $PA$ and $QA$ are conjugate, $QA$
satisfies  $\mathbf{DivSyl}(p)$. Now by Lemma \ref{nupfactor}, $\nu_p(G)=\nu_p(G/A)\cdot \nu_p(QA)$, and
\begin{equation*}
\nu_p(H)=\nu_p(H/(H\cap A))\cdot \nu_p((Q\cap H)(H\cap A))=\nu_p(HA/A)\cdot \nu_p((Q\cap H)(H\cap A)).
\end{equation*}
Since $HA/A\leq G/A$ and $(Q\cap H)(H\cap A)\leq QA$, the conditions of the lemma imply that $\nu_p(HA/A)$ divides
$\nu_p(G/A)$ and $\nu_p((Q\cap H)(H\cap A))$ divides $\nu_p(QA)$. This implies that $\nu_p(H)$ divides~$\nu_p(G)$.
\end{proof}

The following lemma is proven in \cite{Navarro}, but now it can be easily derived from the above two lemmas and we provide
here an alternating proof of it.

\begin{Lemma}\label{psolvable} {\em \cite[Theorem~A]{Navarro}}
Let $H$ be a subgroup of a $p$-solvable group $G$. Then $\nu_p(H)$ divides $\nu_p(G)$.
\end{Lemma}

\begin{proof}
Assume by contradiction that $G$ is a counterexample of minimal order, $H\leq G$ is chosen so that $\nu_p(H)$ does not
divides $\nu_p(G)$ and $A$ is a minimal normal subgroup of $G$. Let $P\in Syl_p(G)$ with $P\cap H\leq Syl_p(H)$. By
Lemma \ref{nupfactor}, we have $$\nu_p(G)=\nu_p(G/A)\cdot \nu_p(PA)$$ and $$\nu_p(H)=\nu_p(HA/A)\cdot \nu_p((H\cap
P)(H\cap A)).$$ By induction, $\nu_p(HA/A)$ divides $\nu_p(G/A)$. Since $G$ is $p$-solvable, $A$ is either a
$p$-group or a $p'$-group.

If $A$ is a $p$-group, then $\nu_p((H\cap P)(H\cap A))=1$, so $\nu_p(H)$
divides $\nu_p(G)$.

If $A$ is a $p'$-group, then $\vert A\vert$ and $\vert P\vert$ are corpime,
in particular $P\cap A=1$. Since $P$ normalizes $A$, we have that
$[N_A(P),P]\leq A\cap P=1$, so $N_A(P)=C_A(P)$. The known properties of
coprime action (for example, see \cite[Excersize~3E.4]{Isaacs}) and Lemma 1
implies that $\nu_p((H\cap P)(H\cap A))=\vert (H\cap A):C_{H\cap A}(H\cap
P)\vert$ divides $\vert A:C_A(H\cap P)\vert$, while $\vert A:C_A(H\cap
P)\vert$ divides $\vert A:C_A(P)\vert=\nu_p(PA)$. Thus $\nu_p(H)$ divides
$\nu_p(G)$, which contradicts the choice of $G$ and $H$.
\end{proof}

Recall that a group $A$ is said to be an almost simple group if there exists a simple nonabelian group $S$ such that
$S\leq A\leq \Aut(S).$ By Schreier conjecture $\Aut(S)/S$ is solvable for each simple group $S$.

\begin{Cor}\label{almsimple}
Let $A$ be an almost simple group with simple socle $S$ and $P\in Syl_p(A)$. If $PS$ satisfies $\mathbf{DivSyl}(p)$,
then $A$ satisfies $\mathbf{DivSyl}(p)$.
\end{Cor}

\begin{proof}
This follows from Lemmas \ref{extension} and \ref{psolvable} since $A/S$ is solvable.
\end{proof}

\begin{Lemma}\label{inducedautomoprhisms} {\em \cite[Theorem~1]{VdoInd}}
Let $$1=G_0<G_1<\ldots<G_n=G$$ be an $(rc)$-series of $G$ and denote the section $G_i/G_{i-1}$ by $S_i$. Let
$$1=H_0<H_1<\ldots<H_n=G$$ be a composition series of $G$. Then there exist a permutation $\sigma\in \Sym_n$ such that
$\Aut_G(H_i/H_{i-1})\leq \Aut_G(S_{i\sigma})$ for $i=1,\ldots,n$. Moreover, if the second series is an $(rc)$-series,
then $\sigma$ can be chosen so that $\Aut_G(H_i/H_{i-1})\simeq\Aut_G(S_{i\sigma})$ for $i=1,\ldots,n$.
\end{Lemma}

The following simple lemma plays technical role.

\begin{Lemma}\label{normpbypdashbyp}
Let $G$ be a finite group, possessing a normal series $$1\lhd H\lhd K\lhd G$$ such that $H$ is a $p$-group, $K/H$ is a
$p'$-group and $G/K$ is a $p$-group. Let $Q$ be a $p$-subgroup of $G$ containing $H$ and $P$ be a Sylow $p$-subgroup of
$G$ containing $Q$. Consider a natural homomorphism $$\overline{\phantom{G}}:G\rightarrow \overline{G}=G/H.$$ Then
$\overline{N_K(Q)}=N_{\overline{K}}(\overline{Q})=C_{\overline{K}}(\overline{Q})$ and $N_K(P)\leq N_K(Q)$.
\end{Lemma}

\begin{proof}
$\overline{N_K(Q)}\subseteq N_{\overline{K}}(\overline{Q})$ is evident.
Consider $\bar{x}\in N_{\overline{K}}(\overline{Q})$. Then every inverse
image $x$ of $\bar{x}$ normalizes the full inverse image of $\overline{Q}$.
Since $H\leq Q$, the full inverse image of $\overline{Q}$ equals $Q$, so $x$
normalizes $Q$. Consequently, $N_{\overline{K}}(\overline{Q})\subseteq
\overline{N_K(Q)}$.  By the conditions of the lemma,
$\overline{K}=O_{p'}(\overline{G})$, so $\overline{K}\unlhd\overline{G}$ and
$\vert \overline{K}\vert$ and $\vert \overline{Q}\vert$ are coprime. Thus
$N_{\overline{K}}(\overline{Q})=C_{\overline{K}}(\overline{Q})$. Since $H\leq
P$ and $P$ is a $p$-group the identity
$\overline{N_K(P)}=N_{\overline{K}}(\overline{P})=C_{\overline{K}}(\overline{P})$
holds. Clearly $C_{\overline{K}}(\overline{P})\leq
C_{\overline{K}}(\overline{Q}).$ It follows that $N_K(P)\leq N_K(Q)$.
\end{proof}

In this paper we often consider subgroups of a permutation wreath product, so we fix notations for such groups. If $L$
is a group and $K$ is a subgroup of the symmetric group $\Sym_k$, then $L\wr K$ denotes the permutation wreath product.
By definition $L\wr K$ is isomorphic to a semidirect product $(L_1\times\ldots\times L_k)\rtimes K$, where $L_1\simeq
\ldots\simeq L_k\simeq L$ and $K$ permutes the $L_i$-s, that is arbitrary $\sigma\in K$ acts by
\begin{equation*}\sigma:(g_1,\ldots,g_k)\mapsto (g_{1\sigma^{-1}},\ldots,g_{k\sigma^{-1}}).\end{equation*}
We always denote $L_1\times\ldots\times L_k$ by $\mathbf{L}$.

Suppose that $G$ is a subgroup of $L\wr K$.  Let $\rho:G\rightarrow \Sym_k$ be the homomorphism corresponding to the action of $G$ on $\{L_1,\ldots,L_k\}$ by
conjugation.  Each element $g$ of $G$ as an element from $L\wr K$ can be uniquely written as
\begin{equation}\label{canonicalformwreath}
g=(g_1,\ldots,g_k)\rho(g),\text{ where }g_i\in L_i,\text{ and }\rho(g)\in\Sym_k\text{ permutes the }g_i\text{-s}.
\end{equation}
Given $g=(g_1,\ldots,g_k)\rho(g)\in L\wr K$ define $\pi_i(g)$ to be equal to $g_i$. In particular, the restriction of
$\pi_i$ on $\mathbf{L}$ is the natural projection on~$L_i$.

Now if $S$ is a normal subgroup of $L$ then we can construct a normal subgroup $\mathbf{S}=S_1\times\ldots\times S_k$ of
$L\wr K$ such that $S_1\simeq\ldots\simeq S_k\simeq S$ and
$(L\wr K)/\mathbf{S}\simeq (L/S)\wr K$. In order to construct $\mathbf{S}$ we need to take one representative $L_j$ from
every $K$-orbit on $\{L_1,\ldots L_k\}$, the in each representative we need to take $S_j\unlhd L_j$ such that
$S_j\simeq S$ and $L_j/S_j\simeq L/S$, and finally define $\mathbf{S}$ to be equal to the normal closure of the subgroup
generated by $S_j$-s. From the construction it is evident, that if $K$ is transitive, then $\mathbf{S}$ is uniquely defined
by the choice of $S_1\leq L_1$.

\begin{Lemma}\label{inclusionwreathproduct}
Assume that $T$ is a unique minimal normal subgroup of $G$. Assume also that $T$ is nonabelian, that is,
$T=S_1\times\ldots\times S_k$, where $S_1,\ldots,S_k$ are  isomorphic (moreover, $G$-conjugate) nonabelian simple
groups. Let $\rho:G\rightarrow \Sym_k$ be the permutation representation corresponding to the action of $G$ on
$\{S_1,\ldots,S_k\}$ by conjugation. Then there exists an injective homomorphism
$$\varphi:G\rightarrow \Aut_G(S_1)\wr \rho(G)=\left(\Aut_G(S_1)\times\ldots\times \Aut_G(S_k)\right)\rtimes \rho(G).$$
Moreover, if $H\leq G$ and $HT=G$, then $\varphi$ can be chosen so that for every $H\leq X\leq G$ we have
$$\varphi(X)\leq \Aut_X(S_1)\wr\rho(X)=\Aut_X(S_1)\wr \rho(G)=\left(\Aut_X(S_1)\times\ldots\times
\Aut_X(S_k)\right)\rtimes \rho(X).$$
\end{Lemma}

\begin{proof}
The existence of $\varphi$ is known (see, for example,
\cite[Lemma~2]{VdoZen}). Hence,  we only need to show that $\varphi(X)\leq
\Aut_X(S_1)\wr \rho(X)$. Here we provide the construction of $\varphi$ and
remain all technical details for the reader.

First notice that $\rho(G)$ is transitive since $T$ is a minimal normal subgroup. By the Fundamental counting
principle (see \cite[Theorem~1.4]{Isaacs}), the action of $G$ on $\{S_1,\ldots,S_k\}$ is equivalent to the action of
$G$ on right cosets of $N_G(S_1)$ by right multiplication. By \cite[Theorem~IV.1.4 a)]{Huppert}, there exists an injective
homomorphism
\begin{equation*}
  \psi:G\rightarrow \left(N_G(S_1)\times\ldots\times N_G(S_k)\right)\rtimes\rho(G)=N_G(S_1)\wr\rho(G),
\end{equation*}
defined in the following way. Fix a right coset representatives  $r_1,\ldots,r_k$ of $N_G(S_1)$ in  $G$. For given
$g\in G$,  define the element $n_i(g)$ of $N_G(S_1)$ by
\begin{equation*}
  r_ig=n_i(g)r_{i\rho(g)}.
\end{equation*}
Then $\psi$ maps $g$ to $(n_1(g),\ldots,n_k(g))\rho(g)$.

Consider a normal subgroup $C_G(S_1)$ of $N_G(S_1)$. Since $\rho(G)$ is transitive, $C_G(S_1)\times\ldots\times C_G(S_k)$
is a normal subgroup. Consider a natural homomorphism
\begin{equation*}
  \theta:N_G(S_1)\wr\rho(G)\rightarrow \left(N_G(S_1)\wr\rho(G)\right)/\left(C_G(S_1)\times\ldots\times C_G(S_k)\right).
\end{equation*}
Since $T$ is the unique minimal normal subgroup of $G$ and $T$ is nonabelian, we have ${\psi(G)\cap Ker(\theta)=1}$, while
$\theta(N_G(S_1)\wr\rho(G))=\Aut_G(S_1)\wr\rho(G)$ by definition.

Now assume that $H\leq G$ is chosen so that $HT=G$. Then the right coset representatives  $r_1,\ldots,r_k$ of $N_G(S_1)$ can be
chosen from $H$. Hence if $H\leq X\leq G$, then for every $x\in X$ and $i=1,\ldots,k$, we have
$n_i(x)=r_ixr_{i\rho(x)}^{-1} \in X\cap N_G(S_1)=N_X(S_1)$. Therefore $\psi(X)\leq N_X(S_1)\wr\rho(X)$, and so
$\varphi(X)=\theta\psi(X)\leq
\Aut_X(S_1)\wr\rho(X)$.
\end{proof}

\begin{Lemma}\label{numpwreath}
Let $L$ be a group with a normal subgroup $S$ such that $L/S$ is a $p$-group. Let $K$ be a transitive $p$-subgroup of
$\Sym_k$. Choose $S_1\leq L_1$ with $S_1\simeq S$ and $L_1/S_1\simeq L/S$, then $\mathbf{S}=S_1^{L\wr K}=S_1\times \ldots\times S_k$.
Let $G$ be a subgroup of
$L\wr K$ satisfying the following conditions:
\begin{itemize}
\item[{\em (a)}] $\mathbf{S}\leq G$;
\item[{\em (b)}] $G\mathbf{L}=L\wr K$;
\item[{\em (c)}] $\pi_i(N_G(S_i))=L_i$.
\item[{\em (d)}] for some $Q\in Syl_p(G)$, if $P_i=\pi_i(N_Q(S_i))$, then
    $Q\leq (P_1\times\ldots\times P_k)\rtimes K$.
\end{itemize}
Then the restriction of $\pi_i$ on $N_G(S_i)$ is a homomorphism and $\nu_p(G)=\vert S\vert_{p'}^{k-1}\cdot\nu_p(L)$.
\end{Lemma}

\begin{proof}
Since $G\mathbf{L}=L\wr K$, we have $\rho(G)=G\mathbf{L}/\mathbf{L}\simeq
G/(G\cap \mathbf{L})\simeq K$, in particular, $G/(G\cap \mathbf{L})$ is a
$p$-group. Hence $G=Q (G\cap \mathbf{L})$ and $L\wr K=Q\mathbf{L}$, in
particular $\rho(Q)=K$. Now $g=(g_1,\ldots,g_k)\rho(g)\in N_{L\wr K}(S_i)$ if
and only if $i\rho(g)=i$. For every $g,h\in L\wr K$ we have $\pi_i(g)\cdot
\pi_i(h)=g_i\cdot h_i$. On the other hand, if we take $g,h\in N_{L\wr
K}(S_i)$, then
\begin{multline*}
  \pi_i(gh)=\pi_i((g_1,\ldots,g_k)\rho(g)\cdot (h_1,\ldots,h_k)\rho(h))\\= \pi_i((g_1,\ldots,g_k)\cdot (h_1,\ldots,h_k)^{\rho(g)^{-1}}\cdot\rho(g) \rho(h))\\=
  \pi_i((g_1,\ldots,g_k)\cdot (h_{1\rho(g)},\ldots,h_i,\ldots,h_{k\rho(g)})\cdot\rho(g) \rho(h))\\= \pi_i((g_1\cdot h_{1\rho(g)},\ldots, g_i\cdot h_i,\ldots g_k
  \cdot h_{k\rho(g)})\cdot\rho(g)\rho(h)))=g_i\cdot h_i,
\end{multline*}
thus the restriction of $\pi_i$ on $N_G(S_i)$ is a homomorphism.

Since $L/S$ is a $p$-group and $K$ is a $p$-group, we have that $(L\wr K)/\mathbf{S}$
is a $p$-group. Hence $G/\mathbf{S}$ is a $p$-group and $Q\mathbf{S}=G$. Denote $\pi_i(Q\cap \mathbf{S})$ by $Q_i$. Clearly $Q\cap S_i\leq Q_i$.
On the other hand, since $S_i$ is subnormal in $G$, $Q\cap S_i$ is a Sylow $p$-subgroup of $S_i$, so $Q\cap S_i=Q_i$.
Since $\mathbf{S}$ normalizes $S_i$ and $G=Q\mathbf{S}$, we have $N_G(S_i)=N_Q(S_i)\mathbf{S}$, therefore $\pi_i(N_G(S_i))=P_iS_i=L_i$. But since
\begin{multline*}
\vert P_i\vert=\vert P_i/(P_i\cap S_i)\vert \cdot \vert P_i\cap S_i\vert=\\
\vert P_i/(P_i\cap S_i)\vert \cdot \vert \pi_i(N_Q(S_i))\cap S_i\vert= \vert P_i/(P_i\cap S_i)\vert \cdot \vert \pi_i(N_Q(S_i)\cap S_i)\vert=\\
\vert P_i/(P_i\cap S_i)\vert \cdot \vert \pi_i(Q\cap S_i)\vert = \vert L_i/S_i\vert\cdot \vert Q_i\vert=\vert L_i\vert_p,
\end{multline*}
we obtain $P_i\in Syl_p(L_i)$. Moreover, since $Q$ as a subgroup of $G$ permutes $S_i$-s, we have that $Q$ permutes
$N_Q(S_i)$-s, and so $Q$ permutes $P_i$-s. This implies that  $Q$ normalizes $P_1\times \ldots\times P_k$. Consequently, $\rho(Q)=K$
normalizes $P_1\times \ldots\times P_k$ and the subgroup $P=(P_1\times \ldots\times P_k)\rtimes K$ in condition (d) of
the lemma is correctly defined. Since $P_i\in Syl_p(L_i)$ and $K$ is a $p$-group, we also have that $P\in Syl_p(L\wr
S)$. Then by condition (d), $Q\leq P$, and so $Q=P\cap G$.

Since $G/S$ is a $p$-group, Lemma \ref{nupfactor} implies that
\begin{equation*}
\nu_p(G)=\vert \mathbf{S}:N_{\mathbf{S}}(Q)\vert
\end{equation*}
so
\begin{equation*}
  \nu_p(G)=\frac{\vert \mathbf{S}\vert_{p'}}{\vert N_{\mathbf{S}}(Q)\vert_{p'}}.
\end{equation*}

Since  $L_1=P_1S_1$, by Lemma \ref{nupfactor} again, we have
\begin{equation}\label{step6identity}
\nu_p(L_1)=\vert S_1:N_{S_1}(P_1)\vert=\frac{\vert S_1\vert_{p'}}{\vert N_{S_1}(P_1)\vert_{p'}}.
\end{equation}
If  we can show that
\begin{equation}\label{principalstep6}
  \vert N_{\mathbf{S}}(Q)\vert_{p'}=\vert N_{S_1}(P_1)\vert_{p'},
\end{equation}
 then \eqref{step6identity} and \eqref{principalstep6} would imply
\begin{equation*}%\label{equation4}
\nu_p(G)=\frac{\vert \mathbf{S}\vert_{p'}}{\vert N_{\mathbf{S}}(Q)\vert_{p'}}=\frac{\vert \mathbf{S}\vert_{p'}}{\vert S_1\vert_{p'}}\cdot
\frac{\vert S_1\vert_{p'}}{\vert N_{S_1}(P_1)\vert_{p'}}=\vert S_1\vert_{p'}^{k-1}\cdot
\nu_p(L_1)
\end{equation*}
and so the lemma follows.

Thus we remain to prove \eqref{principalstep6}.

Denote $N_{S_i}(P_i)$ by $X_i$ and let $\mathbf{X}=X_1\times\ldots\times X_k$. Clearly $Q_i=P_i\cap S_i\unlhd
X_i$, and since $Q_i$ is a Sylow $p$-subgroup of $S_i$, $Q_i$ is a (normal) Sylow $p$-subgroup of $X_i$. Therefore
$Q_1\times\ldots\times Q_k=O_p(\mathbf{X})\in Syl_p(\mathbf{X})$ and $\mathbf{X}$ is an extension of a $p$-group by a $p'$-group. We  show that
$P$ normalizes $\mathbf{X}$ and that $N_{\mathbf{S}}(Q)$ lies in $\mathbf{X}$. By construction,  $P$ permutes $P_i$-s, so $P$ permutes $X_i$-s,
therefore $P$ normalizes $X_1\times\ldots\times X_k=\mathbf{X}$. In particular, $\rho(P)=\rho(Q)=K$ normalizes $\mathbf {X}$ and $\mathbf{X}K=X_1\wr K$.
 In order to show that $N_{\mathbf{S}}(Q)\leq \mathbf{X}$,  consider
$x=(x_1,\ldots,x_k)\in N_{\mathbf{S}}(Q)$. Since $x$ normalizes $Q$ and $S_i$, we obtain that $x$ normalizes $N_Q(S_i)$ for
$i=1,\ldots,k$. Take $y=(y_1,\ldots,y_k)\rho(y)\in N_Q(S_i)$.
Then $i\rho(y)=i$. Recall that by condition (d) of the lemma $Q\leq P$, so $y_i\in P_i$ for $i=1,\ldots,k$. Since $x_j$ centralizes $y_t$ for
$j\not=t$, we obtain
\begin{multline*}
  y^x=(y_1^{x_1},\ldots,y_k^{x_k})\cdot \rho(y)^x= (y_1^{x_1},\ldots,y_k^{x_k})\cdot (x_1^{-1},\ldots,x_k^{-1})\cdot
  (x_1,\ldots,x_k)^{\rho(y)^{-1}}\cdot \rho(y)\\
  =(y_1^{x_1},\ldots,y_k^{x_k})\cdot (x_1^{-1},\ldots,x_i^{-1},\ldots,x_k^{-1})\cdot (x_{1\rho(y)},\ldots,x_i,\ldots, x_{k\rho(y)})\cdot \rho(y)\\
  = (y_1^{x_1}\cdot x_1^{-1}\cdot x_{1\rho(y)},\ldots,y_i^{x_i},\ldots y_k^{x_k}\cdot x_k^{-1}\cdot x_{k\rho(y)})\cdot \rho(y).
\end{multline*}
Thus $\pi_i(y^x)=y_i^{x_i}\in P_i$, i.e. $x_i$ lies in $X_i$ and so $N_{\mathbf{S}}(Q)$ lies in $\mathbf{X}$.

Consider $\mathbf{X}P\leq L\wr K$. Since $Q_1\times\ldots\times Q_k=O_p(\mathbf{X})\in Syl_p(\mathbf{X})$, the group $\mathbf{X}P$ has a normal series
$$1\lhd Q_1\times\ldots\times Q_k=O_p(\mathbf{X})\unlhd \mathbf{X}\unlhd \mathbf{X}P.$$  By Lemma \ref{normpbypdashbyp},
$\overline{N_{\mathbf{X}}(Q)}=N_{\overline{\mathbf{X}}}(\overline{Q})=C_{\overline{\mathbf{X}}}(\overline{Q})$, where
$$\overline{\phantom{G}}:P\mathbf{X}\rightarrow \overline{P\mathbf{X}}=P\mathbf{X}/(Q_1\times\ldots\times Q_k)$$ is a natural homomorphism.

Now consider the corresponding natural homomorphisms
$$\overline{\phantom{G}}:(X_iP_i)\rightarrow \overline{X_iP_i}=(X_iP_i)/Q_i$$ (we simply restrict the above homomorphism
on $X_iP_i$) for $i=1,\ldots,k$. Then each $x=(x_1,\ldots,x_k)\in \mathbf{X}$
is mapped to $\bar{x}=(\bar{x}_1,\ldots,\bar{x}_k)$ and
$y=(y_1,\ldots,y_k)\rho(y)\in Q$ is mapped to
$(\bar{y}_1,\ldots,\bar{y}_k)\rho(y)$. Since $X_i$ normalizes $P_i$, we
obtain that $\overline{X}_i$ centralizes $\overline{P_i}$, i.e.
$(\bar{x}_1,\ldots,\bar{x}_k)$ centralizes $(\bar{y}_1,\ldots,\bar{y}_k)$. So
\begin{multline*}
\bar{x}^{\bar{y}}=(\bar{x}_1,\ldots,\bar{x}_k)^{(\bar{y}_1,\ldots,\bar{y}_k)\rho(y)}=(\bar{x}_1^{\bar{y}_1},\ldots,\bar{x}_k^{\bar{y}_k})^{\rho(y)}=\\
(\bar{x}_1,\ldots,\bar{x}_k)^{\rho(y)}=(\bar{x}_{1\rho(y)^{-1}},\ldots,\bar{x}_{k\rho(y)^{-1}}).
\end{multline*}
Since $\rho(Q)=K\leq\Sym_k$ is transitive, we obtain that
$$C_{\overline{\mathbf{X}}}(\overline{Q})=C_{\overline{\mathbf{X}}}(\rho(Q))=\{(\bar{x}_1,\ldots,\bar{x}_k)\mid
\bar{x}_1=\ldots=\bar{x}_k\}.$$ Thus $$\vert N_{\mathbf{S}}(Q)\vert_{p'}=\vert \overline{N_{\mathbf{X}}(Q)}\vert=C_{\overline{X}}(\overline{Q})=\vert \overline{X_1}\vert
= \vert X_1\vert_{p'}=\vert N_{S_1}(P_1)\vert_{p'}$$ and \eqref{principalstep6} holds.
\end{proof}

\begin{Lemma}\label{numberindirectproduct}
Assume $T$ is a normal subgroup of $G$ such that $T=G_1\times\ldots\times G_k$, and $G/T$ is a $p$-group. Assume also
that $G$ acts on $\{G_1,\ldots,G_k\}$ by conjugation, i.e. for every $x\in G$ and $i=1,\ldots,k$ there exists
$j\in\{1,\ldots,k\}$ with $G_i^x=G_j$. Denote by $\pi_i$ the natural projection $\pi_i:T\rightarrow G_i$. Let $H$ be a
subgroup of $G$ such that for every $i=1,\ldots,k$ we have $\pi_i(H\cap T)$ equals either $G_i$ or $1$, that is, $H\cap T$
is a subdirect product of some of $G_i$-s.
Then $\nu_p(H)$ divides $\nu_p(G)$.
\end{Lemma}

\begin{proof}
Suppose that this lemma is false and consider a counterexample $(G, H)$ with
$|G|$ minimal. Then $H\leq G$ is such that $H\cap T$ is a subdirect product
of some of $G_i$-s and  $\nu_p(H)$ does not divides $\nu_p(G)$.

Let $\Omega_1,\ldots,\Omega_s$ be $G$-orbits of $\{G_1,\ldots, G_k\}$ and
$G_{i_j}$ be a fixed element of $\Omega_j$. In order to simplify notation, we
assume $G_{i_1}=G_1$.  Assume that $L_1$ is a minimal characteristic subgroup
of $G_1$, then $L_1$ is a direct product of isomorphic simple groups, and
either $L_1$ is elementary abelian or $L_1$ is a direct product of isomorphic
nonabelian simple groups. For every $i=2,\ldots,s,$ fix the corresponding
characteristic subgroup $L_{i_j}$ of $G_{i_j}$ so that $L_1\simeq L_{i_j}$
and $G_1/L_1\simeq G_{i_j}/L_{i_j}$. Denote by $L$ the normal closure of
$\langle L_1,\ldots, L_{i_s}\rangle $ in $G$, that is, $L=\langle
L_{i_j}^x\mid x\in G, j=1,\ldots,s\rangle$. Evidently $L$ is a direct product
of some groups isomorphic to $L_1$. In particular, $L$ is abelian if and only
if $L_1$ is abelian and $L$ is a direct product of isomorphic nonabelian
simple groups if and only if $L_1$ is a direct product of isomorphic
nonabelian simple groups. Since $G$ acts transitively by conjugation on each
$G$-orbit, we obtain that $L=\pi_1(L)\times\ldots\times \pi_k(L)$,
$\pi_1(L)\simeq\ldots\simeq\pi_k(L)$, and $G$ also acts by conjugation on
$\{\pi_1(L),\ldots,\pi_k(L)\}$. Hence, if we consider the natural
homomorphism
$$\overline{\phantom{G}}:G\rightarrow \overline{G}=G/L,$$
then
$\overline{T}=\overline{G}_1\times\ldots\times\overline{G}_k=(G_1/\pi_1(L))\times\ldots\times
(G_k/\pi_k(L))$. Clearly $\overline{G}$ acts by conjugation on
$\{\overline{G}_1,\ldots,\overline{G}_k\}$ and
$\overline{G}/\overline{T}\simeq G/T$ is a $p$-group. This shows that
$\overline{G}$ satisfies the conditions of the lemma. We also have
\begin{equation*}
  \pi_i(\overline{H}\cap \overline{T})=\pi_i(\overline{HL\cap T}),
\end{equation*}
so $\pi_i(\overline{H}\cap \overline{T})=\overline{G_i}$ if $\pi(H\cap T)=G_i$ and
$\pi_i(\overline{H}\cap\overline{T})=1$ if $\pi_i(H\cap T)=1$. Thus $\overline{H}$ also satisfies conditions of the
lemma. By induction, $\nu_p(\overline{H})$ divides $\nu_p(\overline{G})$.

Let $Q$ be a Sylow $p$-subgroup of $H$ and $P$ be a Sylow $p$-subgroup of $G$ containing $Q$. By Lemma \ref{nupfactor},
$\nu_p(H)=\nu_p(Q(H\cap L))\cdot \nu_p(H/(H\cap L))=\nu_p(Q(H\cap L))\cdot\nu_p(\overline{H})$ and
$\nu_p(G)=\nu_p(PL)\cdot \nu_p(G/L)=\nu_p(PL)\cdot \nu_p(\overline{G})$. Thus we need to show that $\nu_p(Q(H\cap L))$
divides $\nu_p(PL)$ in order to obtain a contradiction with the choice of $(G, H)$.

Assume first that $L$ is abelian. Then $PL$ is solvable, hence $\nu_p(Q(H\cap L))$ divides $\nu_p(PL)$ by
Lemma~\ref{psolvable}. Thus we may assume that $G_1$ has no characteristic abelian subgroups and so $G_1$ has no
normal abelian subgroups.

Now assume that $L$ is nonabelian. As we note above, $L_1$ is a is a direct product of isomorphic nonabelian simple
groups. So $L$ is also a direct product of isomorphic nonabelian simple groups, say $L=S_1\times\ldots\times S_m$.
Then $G$ acts by conjugation on  $\{S_1,\ldots,S_m\}$. Denote by $\psi_j:L\rightarrow S_j$ the natural projection. If
$\pi_i(H\cap T)=1$, then clearly $\psi_j(H\cap L)=1$ for every $S_j\leq G_i$. Assume that $\pi_i(H\cap T)=G_i$. Since
$H\cap L\unlhd H$, we obtain that $\pi_i(H\cap L)$ is a normal subgroup of $\pi_i(H\cap T)=G_i$. Moreover, $\pi_i(H\cap
L)$ is a normal subgroup of $\pi_i(L)=\prod_{S_j\leq L_i} S_j$, so $\pi_i(H\cap L)$ is a product of some of $S_j$-s.
Therefore, $H\cap L$  is a subdirect product of $S_j$-s. This shows that $(PL, Q(H\cap L))$ satisfies conditions of
the theorem. If $PL\not=G$, then by induction $\nu_p(Q(H\cap L))$ divides $\nu_p(PL)$. Thus we may assume that $PL=G$.

Assume that  $L$ is not a minimal normal subgroup in $PL$. Then there exists $K\unlhd PL$ such that $K< L$. Since
$L=S_1\times\ldots\times S_m$ is a direct product of isomorphic simple groups, we obtain $K=S_{i_1}\times\ldots\times
S_{i_t}$ for some $\{i_1,\ldots,i_t\}\subseteq \{1,\ldots,m\}$. Consider a natural homomorphism $\varphi:PL\rightarrow
PL/K$. Clearly $(\varphi(PL), \varphi(Q(H\cap L)))$ satisfies conditions of the lemma, so by induction
$\nu_p(\varphi(Q(H\cap L)))$ divides $\nu_p(\varphi(PL))$. It is also clear that $(PK, Q(H\cap K))$ also satisfies
the conditions of the lemma, so by induction $\nu_p(Q(H\cap K))$ divides $\nu_p(PK)$. By Lemma \ref{nupfactor},
$\nu_p(Q(H\cap L))=\nu_p(\varphi(Q(H\cap L)))\cdot \nu_p(Q(H\cap K))$ and $\nu_p(PL)=\nu_p(\varphi(PL))\cdot \nu_p(PK)$,
so $\nu_p(Q(H\cap L))$ divides $\nu_p(PL)$. Hence, we may assume that $L$ is a minimal normal
subgroup of~$PL$.

Now we have that $G=PL$, where $P\in Syl_p(G)$, $L=S_1\times\ldots\times S_m$
is a minimal normal subgroup of $G$, and $S_1,\ldots,S_m$ are isomorphic
(even $G$-conjugate) nonabelian simple groups. By Lemma
\ref{inclusionwreathproduct}, there exists an embedding $\varphi:G\rightarrow
\Aut_G(S_1)\wr \rho(G)$, where $\rho:G\rightarrow \Sym_m$ is the permutation
representation corresponding to the action of $G$ on $\{S_1,\ldots,S_m\}$ by
conjugation. Moreover, since $G=PL$, $\varphi$ can be chosen so that for
every $P\leq X\leq G$ we have \begin{equation}\label{incl}\varphi(X)\leq
\Aut_X(S_1)\wr \rho(X)=\Aut_X(S_1)\wr\rho(G).\end{equation} We show that
$\varphi(G)$ as a subgroup of $\Aut_G(S_1)\wr\rho(G)$ satisfies conditions of
Lemma~\ref{numpwreath}.

Since $S_1\times \ldots\times S_m\leq Ker\ \rho$, we obtain that
$\rho(G)=\rho(P)$ is a transitive $p$-subgroup of $\Sym_m$. Now $S_2\times
\ldots \times S_m$ centralizes $S_1$, so $\Aut_G(S_1)/S_1$ is a $p$-group.
Thus $\Aut_G(S_1)$ with normal subgroup $S_1$ plays the role of group $L$ in
the notations of Lemma \ref{numpwreath}, while $\rho(G)$ is a group $K$ in
the notations of Lemma~\ref{numpwreath}.

Since $\varphi(S_1\times\ldots\times S_m)=S_1\times\ldots\times S_m$, we
obtain that $S_1\times\ldots\times S_m\leq \varphi(G)$ and condition (a) of
Lemma \ref{numpwreath} is satisfied. Clearly
$\varphi(G)(\Aut_G(S_1)\times\ldots\times\Aut_G(S_m))=\Aut_G(S_1)\wr\rho(G)$,
so condition (b) of Lemma \ref{numpwreath} is also satisfied. By definition,
$\pi_i(N_G(S_i))=\Aut_G(S_i)$, thus condition (c) is true. Finally, by
\eqref{incl} we have $\varphi(P)\leq
(\Aut_P(S_1)\times\ldots\times\Aut_P(S_m))\rtimes \rho(P)$ and by definition
$\Aut_P(S_i)=\pi_i(N_P(S_i))$, so condition (d) is also true. Thus
\begin{equation}\label{nupG}
\nu_p(G)=\vert S_1\vert_{p'}^{m-1}\cdot \nu_p(\Aut_G(S_1)).
\end{equation}
On the other hand, by Lemma \ref{nupfactor}, $\nu_p(PL)=\vert L:N_L(P)\vert$
and $\nu_p(Q(H\cap L))=\vert (H\cap L):N_{H\cap L}(Q)\vert$.

Assume that $H\cap L\not=L$. Since $L$ is a direct product of nonabelian
simple groups and $H\cap L$ is a subdirect product of some simple factors of
$L$, we obtain that $\vert H\cap L\vert =\vert S_1\vert^s$ for some $s<m$.
Clearly $\vert (H\cap L):N_{H\cap L}(Q)\vert$ divides $\vert
S_1\vert_{p'}^s$, since $\vert S_1\vert_{p'}^s$ is the index of a Sylow
$p$-subgroup of $H\cap L$. Since $s\leq m-1$, \eqref{nupG} implies $\nu_p(H)$
divides $\nu_p(G)$.

Assume finally that  $H\cap L=L$. Then from $Q\leq P$ we obtain $Q\cap
L=P\cap L\in Syl_p(L)$. Hence the group $PN_L(P\cap L)$ has the following
normal series:$$1\lhd (P\cap L)\lhd N_L(P\cap L)\lhd PN_L(P\cap L),$$ where
$P\cap L$ is a $p$-group, $N_L(P\cap L)/(P\cap L)$ is a $p'$-group, and
$PN_L(P\cap L)/N_L(P\cap L)\simeq P/(P\cap N_L(P\cap L))$ is a $p$-group. By
Lemma \ref{normpbypdashbyp}, $N_{N_L(P\cap L)}(Q)\geq N_{N_L(P\cap L)}(P)$.
But, clearly, $N_{H\cap L}(Q)=N_L(Q)$ normalizes $Q\cap L$, so lies in
$N_L(P\cap L)$. Hence
$$N_{H\cap L}(Q)=N_L(Q)=N_{N_L(P\cap L)}(Q)\geq N_{N_L(P\cap L)}(P)= N_L(P),$$
and so $\nu_p(Q(H\cap L))$ divides $\nu_p(PL)$.
\end{proof}

\section{Proof of the main theorem}

Assume this theorem is false and let $G$ be a counterexample of minimal order. We proceed in a series
of steps to get a contradiction.
\medskip

{\slshape Step 1.} If $1\not=A\unlhd G$, then $G/A$ satisfies $\mathbf{DivSyl}(p)$.
\medskip

Indeed, by Lemma \ref{inducedautomoprhisms}, the group of induced automorphisms does not depend on the choice of an
$(rc)$-series. Hence we may assume that $$1=G_0<G_1<\ldots<G_n=G$$ goes through $A$, i.e. $A=G_i$ for some $i$. Consider a
natural homomorphism $\overline{\phantom{G}}:G\rightarrow \overline{G}=G/A$. Clearly
$$\overline{1}=\overline{G}_i<\overline{G}_{i+1}<\ldots<\overline{G}_n=\overline{G}$$ is an $(rc)$-series of
$\overline{G}$. By \cite[Lemma~1.2]{VdCart}, $\Aut_{\overline{G}}(\overline{G}_j/\overline{G}_{j-1})\simeq
\Aut_G(G_j/G_{j-1})$ for $j=i+1,\ldots,n$, so $\overline{G}$ satisfies $\mathbf{DivSyl}(p)$ by induction.
\medskip

{\slshape Step 2.} The solvable radical $S(G)$ (that is the largest solvable
normal subgroup of $G$) of $G$ is trivial.
\medskip

Otherwise $G/S(G)$ satisfies  $\mathbf{DivSyl}(p)$ by Step 1. Moreover, for given $P\in Syl_p(G)$, the group $PS(G)$ satisfies
$\mathbf{DivSyl}(p)$ by Lemma \ref{psolvable}. Hence $G$ satisfies $\mathbf{DivSyl}(p)$ by Lemma \ref{extension}, a
contradiction.
\medskip

{\slshape Step 3.} Let $T=S_1\times\ldots\times S_k$ be a minimal normal subgroup of $G$. Then $G/T$ is a $p$-group.
Consequently, $T$ is the unique minimal normal subgroup of $G$, and  $S_i$ is a nonabelian simple group for $i=1, \ldots, k.$
\medskip

Indeed, let $P$ be a Sylow $p$-subgroup of $G$. If $G/T$ is not a $p$-group, then $PT$ is a proper subgroup of $G$.
By Lemma \ref{inducedautomoprhisms}, we may assume that  $$1=G_0<G_1<\ldots<G_n=G$$ goes through $T$, that is,
$T=G_k$ for some $k\in \{1, \ldots, n-1\}.$  Denote the image of $P\cap N_G(G_i/G_{i-1})$ in $\Aut_G(G_i/G_{i-1})$  by $P_i$ for $i=1,\ldots,k$. By
definition, $P_i(G_i/G_{i-1})=\Aut_{PT}(G_i/G_{i-1})\leq \Aut_G(G_i/G_{i-1})$, and the conditions of the theorem says that
for every $p$-subgroup $Q$ of $P_i(G_i/G_{i-1})$ the group $Q(G_i/G_{i-1})$ satisfies $\mathbf{DivSyl}(p)$. Hence $PT$
satisfies the conditions of the theorem, so by induction, $PT$ satisfies $\mathbf{DivSyl}(p)$. Then Step 1 and Lemma
\ref{extension} imply that $G$ satisfies $\mathbf{DivSyl}(p)$, a contradiction. \medskip

Therefore $G/T$ is a $p$-group and so $G=PT$ for some (hence for every) Sylow
$p$-subgroup $P$ of $G$. Consequently, $G$ has the unique minimal normal subgroup $T=L_1\times\ldots\times L_k$, where all $L_i$-s are
isomorphic nonabelian simple groups by Step 2.
\medskip

{\slshape Step 4.} If $H$ is a subgroup of $G$ containing $T$, then $\nu_p(H)$ divides $\nu_p(G)$.
\medskip

Assume that $T\leq H$. Let $Q$ be a Sylow $p$-subgroup of $H$. Since $G/T$ is
a $p$-group, so is $H/T$,  in particular, $H=QT$. Let $P\in Syl_p(G)$ be
chosen so that $P\cap H=Q$. Then  $P\cap T=Q\cap T\in Syl_p(T)$ is a normal
subgroup of $P$ and of $Q$. Denote $N_T(P\cap T)$ by $K$. Since $P\cap T\in
Syl_p(T)$, we have $P\cap T\in Syl_p(K)$, so $K/(P\cap T)$ is a $p'$-group,
and $PK$ satisfies conditions of Lemma \ref{normpbypdashbyp}. Clearly both
$N_T(Q)$ and $N_T(P)$ lie in $K$, so $N_T(Q)=N_K(Q)$ and $N_T(P)=N_K(P)$. By
Lemma \ref{normpbypdashbyp}, $N_T(Q)=N_K(Q)\geq N_K(P)=N_T(P)$. Finally by
Lemma 1, $\nu_p(H)=\nu_p(QT)=\vert T:N_T(Q)\vert$ and
$\nu_p(G)=\nu_p(PT)=\vert T:N_T(P)\vert$. Thus $\nu_p(H)$ divides~$\nu_p(G).$
\medskip

Denote $\Aut_G(S_i)$ by $L_i$. Since $T$ is a unique minimal normal subgroup
and $T$ is nonabelian, by Lemma \ref{inclusionwreathproduct}, there exists an
injective homomorphism $\varphi:G\rightarrow L_1\wr \rho(G)$, where
$\rho:G\rightarrow \Sym_k$ is the  permutation representation corresponding
to the action of $G$ on $\{S_1,\ldots,S_k\}$ by conjugation. Moreover, since
$S_i$ is a unique minimal normal subgroup of $L_i$ and $\rho(G)$ is
transitive, we have $\mathbf{S}=S_1\times\ldots\times S_k$ is a unique
minimal normal subgroup of $L\wr \rho(G)$. In particular, for every such
inclusion $\varphi$ the identity $\varphi(T)=\mathbf{S}$ holds. Below we
identify $G$ with its image $\varphi(G)\leq L_1\wr \rho(G)$ and $T$ with
$\mathbf{S}$. For every $H\leq G$ denote $\pi_i(H\cap T)$ by $H_i$, where
$\pi_i$ is defined before Lemma \ref{inclusionwreathproduct}. Clearly $H\cap
T$ is a subdirect product of $H_1,\ldots,H_k$.
\medskip

{\slshape Step 5.} If $H\leq G$ is chosen so that $HT=G$, then $\vert H_1\vert=\ldots=\vert H_k\vert$ and $H$
normalizes $H_1\times \ldots\times H_k$.
\medskip

Since $T$ is a minimal normal subgroup of $G$ and $T$ is nonabelian, $G$ acts transitively by conjugation on
$\{S_1,\ldots,S_k\}$.  Since $HT=G$ and $T\leq Ker(\rho)$, we obtain $\rho(H)=\rho(G)$, in particular, for
every $i=1,\ldots,k$, there exists $x_i\in H$ such that $S_1^{x_i}=S_i$.

Let $h\in H\cap T$. Then $h=h_1\cdot\ldots\cdot h_k$, where $h_i=\pi_i(h)\in L_i$. Since $H\cap T\unlhd H$, it
follows that $h^{x_i}\in H\cap T$ for $i=1,\ldots,k$. On the other hand, $h^{x_i}=h_1^{x_i}\cdot\ldots\cdot h_k^{x_i}$
and this identity with the fact that $x_i$ permutes $S_i$-s imply $\pi_i(h^{x_i})=h_1^{x_i}$. Whence
$H_1^{x_i}\subseteq H_i$. The same arguments show that $H_i^{x_i^{-1}}\subseteq H_1$. Hence $H_1^{x_i}=H_i$ and $\vert
H_1\vert=\vert H_i\vert$ for $i=1,\ldots,k$. These arguments also show that arbitrary $x\in H$ permutes $H_i$-s, so $x$
normalizes $H_1\times\ldots\times H_k$.\medskip

{\slshape Step 6.} Assume $H\leq G$ is chosen so that $HT=G$. If  $H\cap T=H_1\times \ldots\times H_k$, then
$\nu_p(H)=\nu_p(\Aut_H(S_1))\cdot\vert H_1\vert_{p'}^{k-1}$ divides $\nu_p(G)$. If  $H\cap T\not =H_1\times \ldots\times H_k$, then
$\nu_p(H)$ divides $\nu_p(\Aut_H(S_1))\cdot\vert H_1\vert_{p'}^{k-1}$, and so $\nu_p(H)$ also divides $\nu_p(G)$.
\medskip

Let $Q$ be a Sylow $p$-subgroup of $H$ and $P\in Syl_p(G)$ is chosen so that $P\cap H=Q$. By Step 3, $G/T$ is a $p$-group. So $H=Q(H\cap
T)$ and $G=QT$. By Lemma \ref{inclusionwreathproduct}, there exists an embedding $\varphi:G\rightarrow \Aut_G(S_1)\wr
\rho(G)$ such that for every $Q\leq X\leq G$ we have $\varphi(X)\leq \Aut_X(S_1)\wr \rho(Q)$, so we choose the embedding $\varphi$ with this property.
Clearly $G$ as a subgroup of $\Aut_G(S_1)\wr\rho(G)$ with Sylow $p$-subgroup $P\leq \Aut_P(S_1)\wr \rho(G)$ satisfies conditions of Lemma
\ref{numpwreath}, so $\nu_p(G)=\vert S_1\vert_{p'}^{k-1}\cdot \nu_p(\Aut_G(S_1))$. If $H\cap T=H_1\times\ldots\times
H_k$, then $H$ as a subgroup of $\Aut_H(S_1)\wr \rho(H)$ with Sylow $p$-subgroup $Q\leq \Aut_Q(S_1)\wr\rho(G)$ satisfies conditions of Lemma \ref{numpwreath}, so
$\nu_p(H)=\vert H_1\vert_{p'}^{k-1}\cdot\nu_p(\Aut_H(S_1))$. Clearly $\vert \Aut_H(S_1)\vert$ divides
$\vert\Aut_G(S_1)\vert$, while $\nu_p(\Aut_H(S_1))$ divides $\nu_p(\Aut_G(S_1))$ by the conditions of the theorem.
Hence $\nu_p(H)$ divides $\nu_p(G)$.

If $H\cap T\not = H_1\times\ldots\times H_k$, then $H$ as a subgroup of $H(H_1\times\ldots\times H_k)$ satisfies
conditions of Lemma \ref{numberindirectproduct}. So $\nu_p(H)$ divides $\nu_p(H(H_1\times\ldots\times H_k))$, and
$\nu_p(H(H_1\times\ldots\times H_k))=\vert H_1\vert_{p'}^{k-1}\cdot\nu_p(\Aut_H(S_1))$ divides $\nu_p(G)$ in
view of the previous paragraph. Hence $\nu_p(H)$ divides $\nu_p(G)$.

\medskip

{\slshape Step 7.} If  $H\leq G$ and $HT\neq G$, then $\nu_p(H)$ divides $\nu_p(G)$.
\medskip

Since nonabelian composition factors of $HT$ and $G$ are the same, like in the proof of Step 3,
we can show that $HT$ satisfies conditions of the theorem, so $HT$ satisfies  $\mathbf{DivSyl}(p)$ by induction.
Therefore $\nu_p(H)$ divides $\nu_p(HT)$. But by Step 4, $\nu_p(HT)$ divides $\nu_p(G)$. Consequently, $\nu_p(H)$ divides $\nu_p(G)$.
\medskip

In view of Steps 4, 6 and 7, we obtain that in any case $\nu_p(H)$ divides
$\nu_p(G)$, which contradicts the choice of $G$. This completes the proof.

\section{Remarks and questions}

In this section we discuss natural question arising after the main theorem. We now can prove the following

\begin{Prop}\label{psolvalmsimple}
Let $G$ be an almost simple group with simple socle $S$ such that $p$ does not divides $\vert S\vert$. Then $G$
satisfies $\mathbf{DivSyl}(p)$.
\end{Prop}

\begin{proof}
By Corollary \ref{almsimple}, we only need to verify that $PS$ satisfies  $\mathbf{DivSyl}(p)$ for $P\in Syl_p(G)$. Let $H$
be a subgroup of $PS$. We may assume that $P\cap H\in Syl_p(H)$. Since $p$ does not divides $\vert  S\vert$, we have
$PS=S\rtimes P$, and so $H=(H\cap S)\rtimes (H\cap P)$.  The known properties of coprime action (see
\cite[Excersize~3E.4]{Isaacs}, for example) implies that $\nu_p((H\cap S)\rtimes (H\cap P))=\vert (H\cap S):C_{H\cap
S}(H\cap P)\vert$ divides $\vert S:C_S(H\cap P)\vert$, while $\vert S:C_S(H\cap P)\vert$ divides $\vert
S:C_S(P)\vert=\nu_p(PS)$. Therefore $\nu_p(H)$ divides $\nu_p(PS)$.
\end{proof}

From Proposition \ref{psolvalmsimple} we immediately get  the following remark

\begin{Remark} The main theorem is the direct generalization of Navarro theorem since all nonabelian
composition factors of a $p$-solvable group are $p'$-groups.
\end{Remark}

\begin{Remark} In conditions of the main theorem, we assume that
$P(G_i/G_{i-1})$ satisfies  $\mathbf{DivSyl}(p)$ for every $p$-subgroup $P$ of $\Aut_G(G_i/G_{i-1})$. Corollary
\ref{almsimple} implies that under this condition every group $L$ with $G_i/G_{i-1}\leq L\leq \Aut_G(G_i/G_{i-1})$
satisfies $\mathbf{DivSyl}(p)$. However the authors do not know any example of almost simple group $L$ with simple
socle $S$ such that $L$ satisfies $\mathbf{DivSyl}(p)$, while $S$ does not. So it is natural to assert

\begin{Conj}\label{simpleDivSyl}
Let $S$ be a simple group satisfying $\mathbf{DivSyl}(p)$. Then every $L$ such that $S\leq L\leq \Aut(S)$ satisfies
$\mathbf{DivSyl}(p)$.
\end{Conj}
\end{Remark}

\begin{Remark}\label{examples} The of the normalizers of Sylow $2$-subgroups in finite simple groups (see \cite{Kond}) implies that in (infinitely) many
almost simple groups the Sylow $2$-subgroups are self-normalizing, so these groups satisfy $\mathbf{DivSyl}(2)$. Also the
classification of subgroups in $\SL_2(p^t)$ implies that $\PSL_2(p^t)$ satisfies $\mathbf{DivSyl}(p)$. Of course, there
are over examples. We provide these examples  just to show that there are may almost simple groups satisfying $\mathbf{DivSyl}(p)$ that are not
$p'$-groups.
\end{Remark}

In connection with Conjecture \ref{simpleDivSyl} and Remark \ref{examples}, the following question arising naturally:

\begin{Prob}
How often almost simple groups satisfy $\mathbf{DivSyl}(p)$?
\end{Prob}

\begin{Prob}
What almost simple groups satisfy $\mathbf{DivSyl}(p)$ for every prime~$p$?
\end{Prob}

\begin{Remark} Another possible questions arise: If we try to consider $\pi$-Hall subgroups instead of Sylow $p$-subgroups. In
\cite{Turull} A. Turull prove that if $G$ is $\pi$-separable, then $\nu_\pi(H)$ divides $\nu_\pi(G)$ for every $H\leq
G$. It is not hard to see that Lemmas \ref{nupfactor}--\ref{psolvable} hold if we consider $\pi$-separable groups and
$\pi$-Hall subgroups, so an alternative proof of Turull theorem can be also obtained by the same arguments. But in
order to extend Turull theorem, one needs to introduce appropriate class of finite groups so that each proper subgroup
satisfies $E_\pi$ at least. We do not go into details in this case, we only mention that such possibility exists and
could be considered later.
\end{Remark}

\end{document}